\documentclass[12pt]{amsart}
\usepackage{amsfonts}
\usepackage{amssymb}

\usepackage{amsmath}
\usepackage{epsf}
\usepackage{epsfig}
\usepackage{colordvi}
\usepackage{fancybox}

\usepackage{colordvi}

\newcommand{\R}{\mathbb R}

\newcommand{\E}{\mathbb E}

\newtheorem{thm}{Theorem}[section]

\newtheorem{prop}[thm]{Proposition}
\theoremstyle{definition}

\newtheorem{example}{Example}[section]

\theoremstyle{remark}

\newcommand{\ds}{\displaystyle}

\begin{document}

\title[ Weierstrass  Representations  of  Lorentzian Minimal Surfaces in $\mathbb R^4_2$]
{Weierstrass Representations  \\ of Lorentzian Minimal Surfaces in $\mathbb R^4_2$}%

\author{Ognian Kassabov, Velichka Milousheva}

\address{Department of Mathematics and Informatics, Todor Kableshkov University of Transport,
158 Geo Milev Str., 1574 Sofia, Bulgaria}
\email{okassabov@abv.bg}

\address{Institute of Mathematics and Informatics, Bulgarian Academy of Sciences,
Acad. G. Bonchev Str. bl. 8, 1113, Sofia, Bulgaria}
\email{vmil@math.bas.bg}


\begin{abstract}

The minimal Lorentzian surfaces in $\mathbb{R}^4_2$ whose first normal space is two-dimensional and whose Gauss curvature $K$ and normal curvature $\varkappa$ satisfy  $K^2-\varkappa^2 >0$ are called minimal Lorentzian surfaces of general type. These surfaces admit canonical parameters and with respect to such parameters  are  determined uniquely up to a motion in $\mathbb{R}^4_2$ by the curvatures $K$ and $\varkappa$ satisfying a system of two natural PDEs.

In the present paper we study minimal Lorentzian surfaces  in $\mathbb{R}^4_2$ and find a Weierstrass representation with respect to isothermal parameters of any minimal surface with two-dimensional first normal space.  We also obtain  a Weierstrass representation with respect to canonical parameters of any minimal Lorentzian surface of general type and solve explicitly the system of natural PDEs expressing any solution to this system by means of four real functions of one variable.

\end{abstract}

\maketitle

\section{Introduction} \label{S:Int}

The study of minimal surfaces is one of the main topics in classical differential geometry. In  the last years, great attention is paid to Lorentzian surfaces in pseudo-Euclidean spaces, since pseudo-Riemannian geometry has many important applications in Physics. 
Minimal Lorentzian surfaces in $\mathbb{C}^2_1$ have been classified by  B.-Y. Chen \cite{Chen2009}.
Classification results for minimal Lorentzian surfaces in  pseudo-Euclidean space $\E^m_s$ with arbitrary dimension $m$ and arbitrary index $s$  are obtained in \cite{Chen2011}. Minimal Lorentzian surfaces in the pseudo-Euclidean 4-space  $\R^4_2$ whose first normal space is two-dimensional and whose Gauss curvature $K$ and normal curvature $\varkappa$ satisfy the inequality $K^2-\varkappa^2 >0$ are studied in \cite{A-M} under the name minimal Lorentzian surfaces of general type. This class of surfaces is characterized  in terms of a pair of smooth functions satisfying a system of two natural partial differential equations. The approach to the study of  minimal Lorentzian surfaces of general type in $\E^4_2$ is based on the introducing of special geometric parameters which are called canonical parameters. 

 A representation of a minimal Lorentzian surface was given by M.P. Dussan and M. Magid in \cite{D-Mag} where they solved the Bj\"orling problem for timelike surfaces in $\R^4_2$ constructing  a special  normal frame and a split-complex representation formula. The Bj\"orling problem for timelike surfaces in the Lorentz-Minkowski spaces $\R^3_1$ and $\R^4_1$ is solved in \cite{Ch-D-Mag} and  \cite{D-Fil-Mag}, respectively. Spinor representation of Lorentzian surfaces in the pseudo-Euclidean 4-space with neutral metric is given in \cite{Bay-Patty}.
In \cite{Patty}, V. Patty gave a generalized Weierstrass representation of a minimal Lorentzian surface in $\R^4_2$ using spinors and Lorentz numbers (also known as para-complex, split-complex, double or hyperbolic numbers)  thus extending the Weierstrass representation of a minimal surface in $\R^3_1$ given by J. Konderak \cite{Kond}.

A special Weierstrass representation, called canonical Weierstrass representation, for a maximal spacelike surface in $\R^4_2$ is obtained in \cite{G-K-2019}. These Weierstrass formulas give explicitly the solutions to the system of PDEs describing maximal spacelike surfaces in $\R^4_2$ in terms of two holomorphic functions. The explicit solving of the system of PDEs and a relation between maximal spacelike surfaces in $\R^4_2$ and maximal spacelike surfaces in $\R^3_1$ is given in \cite{G-K-CR2019} - a paper that is a follow-up to a series of articles on explicit solving of the system of natural PDEs of minimal surfaces in $\R^4$ and minimal spacelike surfaces in $\R^4_1$ \cite{G-K-CR2014, G-K-CR2017}.

\vskip 1mm

In the present paper we study minimal Lorentzian surfaces in $\R^4_2$. In Section \ref{S:W-type formula} we find a Weierstrass representation with respect to isothermal parameters (Theorem \ref{T:3.3}) of any minimal Lorentzian surface with two-dimensional first normal space.   The Weierstrass representation
formula given in Theorem \ref{T:3.3} is equivalent to the representation formula in \cite{D-Mag}. 
But the formula in \cite{D-Mag} uses Lorentz numbers (split-complex numbers), while our formula gives a representation in terms of real functions of one real variable which is much more convenient. Moreover, this Weierstrass 
formula gives a clear analytic form of the surface. 

Further, we obtain a Weierstrass representation
formula with respect to canonical parameters of any minimal Lorentzian surface of general type which describes locally all these surfaces in
terms of four real functions (Theorem \ref{T:4.3}).
Using the canonical Weierstrass representation we solve explicitly the system of 
PDEs that describes minimal Lorentzian surfaces of general type in $\R^4_2$, expressing any solution to this system by means of four real functions of one variable (Theorem \ref{T:4.4}).

In the last section, using the canonical Weierstrass representation formula we give examples of minimal Lorentzian surfaces of general type parametrized by canonical parameters.


\vspace{1. cm}

\setcounter{equation}{0}
\section{Preliminaries} \label{S:Pre}

Let $\R^4_2$ be the four-dimensional pseudo-Euclidean space with the standard flat
metric $\langle \, , \, \rangle$ of signature $(2,2)$, also known as pseudo-Euclidean space with neutral metric, and consider a minimal Lorentzian surface in it. A local parametrization of such 
surfaces is given by the following fundamental result:

\vspace{0.2 cm}
{\bf Theorem A \cite{Chen2009}.}  
{\it Let $S$ be a minimal Lorentzian surface in $\R^4_2$.
Then it can be locally parametrized in the form
$$
	\Psi(u,v)=\beta(u)+\theta(v),
$$
where $\beta$ and $\theta$ are null curves in $\R^4_2$, such that
$\langle\beta'(u),\theta'(v)\rangle \neq 0$ for every $u,\, v$.}

\vspace{0.2 cm}
Concerning generalizations of this theorem for more general situations see \cite{Anciaux},  p. 53 and  \cite{Chen2011}.

The Gauss curvature $K$  and the curvature of the normal connection (the normal curvature) $\varkappa$ of any minimal surface in the Euclidean space $\R^4$ and any spacelike maximal surface in $\R^4_2$ satisfy the inequality  $K^2-\varkappa^2 \geq 0$. The equality case characterizes the so-called  {\it super-conformal} surfaces.
Minimal Lorentz surfaces in the pseudo-Euclidean space $\R^4_2$ can  be divided into three  basic classes:
\begin{itemize}
\item surfaces characterized by $K^2 - \varkappa^2 > 0$; 
\item  surfaces characterized by  $K^2-\varkappa^2 = 0$;
\item surfaces characterized by  $K^2-\varkappa^2 < 0$.
\end{itemize}
 In the case  $K^2-\varkappa^2 > 0$, Lorenzian minimal surfaces in $\R^4_2$ have similar properties to minimal surfaces in $\R^4$ or spacelike maximal surfaces in $\R^4_2$ \cite{Sak}.
Minimal Lorentzian  surfaces  from the first class, i.e. satisfying the inequality  $K^2 - \varkappa^2 > 0$ are studied in \cite{A-M}. In the special case when the first normal space is one-dimensional at each point, the minimal surface is either a degenerate hyperplane, a flat umbilic surface, a quasi-umbilic surface, or lies in a non-degenerate hyperplane. In the case when the first normal space is two-dimensional, the surfaces  are called in \cite{A-M} \textit{minimal  Lorentzian surfaces of general type}, and we will follow this terminology.  Note that minimal  surfaces  satisfying the inequality  $K^2 - \varkappa^2 < 0$ do not exist neither in the Euclidean space $\R^4$ nor in the family of spacelike maximal surfaces in $\R^4_2$ and they need a different approach to be studied with.

\vskip 3mm
In what follows we consider minimal Lorentzian surface of general type, i.e. at any point  $K^2-\varkappa^2 > 0$ 
and the first normal space is two-dimensional.

Let the minimal Lorentzian surface of general type 
$$
	S:  \Psi=\Psi(u,v), \; (u, v) \in \mathcal{D}, \,  \mathcal{D} \subset \R^2
$$
be parametrized by isothermal parameters, i.e.
$E=-G$, $F=0$. Then the coefficients of the first fundamental form are expressed as
$E = f^2(u, v), \; F = 0, \; G = -f^2(u, v)$ for some positive function $f(u, v)$. 
Denote by $x$, $y$ the unit tangent vector fields  in the
directions of $\Psi_u$, resp. $\Psi_v$, i.e. $x=\ds{\frac{\Psi_u}{f}}$, $y=\ds{\frac{\Psi_v}{f}}$. 
Obviously, $\langle x, x \rangle = 1$, $\langle y, y \rangle = -1$, $\langle x, y \rangle = 0$. 
We choose a local  normal frame field $\{e_1, e_2 \}$ such that $\langle e_1, e_1 \rangle = 1$, $\langle e_2, e_2 \rangle = -1$, $\langle e_1, e_2 \rangle = 0$ and define functions $a(u,v), b(u,v), c(u,v), d(u,v)$ by
$$
	\begin{array}{l}
	\vspace{2mm} 
		\sigma(x,x)=a e_1+b e_2; \\
		\vspace{2mm} 
		\sigma(x,y)=c e_1+d e_2; \\
		\vspace{2mm} 
		\sigma(y,y)=a e_1+b e_2,
	\end{array}
$$
where $\sigma$ is the second fundamental tensor. Let $A_1$ and $A_2$ be the shape operators corresponding to $e_1$ and $e_2$, respectively. Then 
$$
\begin{array}{ll}
\vspace{2mm}A_{1}(x)= a x - c y; \qquad & A_{2}(x)= - b x + d y;\\
\vspace{2mm}A_{1}(y)= c x - ay;\qquad   & A_{2}(y)= - d x + b y.
\end{array}
$$
Using the above formulas and the Ricci equation we obtain that the Gauss curvature $K$ and the curvature of the normal connection $\varkappa$ are expressed by
\begin{equation} \label{eq:2.1}
	K= b^2-a^2+c^2-d^2, 
\end{equation}
\begin{equation} \label{eq:2.2}
	\varkappa=2(bc-ad).
\end{equation}

In \cite{A-M} it is proved that in a neighbourhood of any point of a minimal Lorentzian surface of
general type we can introduce a special orthonormal frame field 
$\{x, y, n_1, n_2\}$  such that the vector fields $\sigma(x, x)$ and $\sigma(x, y)$ are collinear to $n_1$ and $n_2$, respectively, i.e.
\begin {equation*} \label{eq:2.3}
\begin {array} {l}
\vspace{2mm}
\sigma(x, x) = \nu n_1;\\
\vspace{2mm}
\sigma(x, y) = \qquad  \mu n_2;\\
\vspace{2mm}
\sigma(y, y) = \nu n_1,
\end {array}
\end {equation*}
where $\langle x, x \rangle = 1$, $\langle y, y  \rangle= -1$, 
$\langle n_1, n_1 \rangle = \varepsilon$, $\langle n_2, n_2 \rangle = -\varepsilon$ 
($\varepsilon = \pm 1$) and $\nu\mu$ never vanishes. The tangent directions determined by the tangent vector fields $x$ and $y$ 
are called \textit{canonical directions} of the surface. 

It follows from 
(\ref{eq:2.1}) and (\ref{eq:2.2})  that the following equality holds:
\begin{equation} \label{eq:2.6}
	K^2-\varkappa^2 = (\mu^2-\nu^2)^2.
\end{equation}
Since we are  studying  surfaces for which $K^2 - \varkappa^2 >0$, we assume that $\mu ^2 - \nu^2 \neq 0$.

Moreover, if 
\begin{equation} \label{eq:2.4}
	E^2=G^2=\frac1{|\mu^2-\nu^2|},  
\end{equation}
the parameters of the surface are called {\it canonical} \cite{A-M}. We shall use the following

\vspace{0.2 cm}
{\bf Theorem B \cite{A-M}.}
{\it The Gauss curvature $K$ and the curvature of the normal connection $\varkappa$ 
of any minimal Lorentzian surface of general type in $\R^4_2$ 
 in canonical parameters satisfy the system
\begin{equation} \label{eq:2.5}
	\begin{array}{l}
	\vspace{2mm}
		\ds\root 4 \of {K^2-\varkappa^2}\Delta^h\ln(K^2-\varkappa^2)=8K  \vspace{2mm}\\ 
		\vspace{2mm}
		\ds\root 4 \of {K^2-\varkappa^2}\Delta^h\ln\frac{K+\varepsilon\varkappa}{K-\varepsilon\varkappa}=4\varepsilon\varkappa
	\end{array}
\end{equation}
where $\Delta^h$ is the hyperbolic Laplacian, $\varepsilon = \pm 1$. Conversely, for any solution $(K(u,v),\varkappa(u,v))$ to this 
system with non-vanishing functions $K$ and $\varkappa$, satisfying $K^2-\varkappa^2>0$,
there exists a unique minimal Lorentzian surface of general type with Gauss curvature $K$ and 
curvature of the normal connection $\varkappa$. 
Moreover, the parameters $(u,v)$ are canonical.}


\vspace{.5 cm}

\setcounter{equation}{0}

\section{The Weierstrass-type formula} \label{S:W-type formula}

The following statement is a trivial consequence of Theorem A:

\vspace{0.2 cm}
\begin{prop}\label{T:3.1}   
Let $\beta$ and $\theta$ be null curves in $\R^4_2$, such that
$\langle\beta',\theta'\rangle \neq 0$.
Then  
\begin{equation*} \label{eq:MinSurf}
	S:  \Psi(u,v)=\beta(u+v)+\theta(u-v)
\end{equation*}
is a minimal Lorentzian surface in isothermal parameters. Conversely, every
minimal Lorentzian surface in $\R^4_2$ can be parametrized locally in this way.
\end{prop}

Now let 
$$
	\beta(u)=\big( \beta_1(u),  \beta_2(u), \beta_3(u), \beta_4(u) \big), \quad u \in J \subset \R
$$
be a null curve in  $\R^4_2$, i.e.
\begin{equation} \label{eq:3.1}
	(\beta_1')^2+(\beta_2')^2-(\beta_3')^2-(\beta_4')^2=0. 
\end{equation}

First we suppose that $(\beta_1')^2=(\beta_3')^2=(\beta_4')^2$ in an interval $J_0 \subset J$. Then $(\beta_2')^2=(\beta_1')^2 $ and hence, there exists a function $f(u), \; u \in J_0$ such that 
$$
	\beta'(u)=f(u)\big( 1,\varepsilon_1,\varepsilon_2,\varepsilon_3\big), \quad \varepsilon_i=\pm 1, \; i = 1, 2, 3.
$$

Now assume  that $(\beta_1')^2-(\beta_3')^2 \neq 0$ at a point $u_0 \in J$. Then the function
$(\beta_1')^2-(\beta_3')^2$ doesn't vanish in an interval $J_0 \subset J$, $u_0 \in J_0$.
We introduce in $J_0$ the following functions:
$$
	f=\frac{\beta'_1+\beta'_3}{2}, \qquad 	
	g=\frac{\beta'_2+\beta'_4}{\beta'_1+\beta'_3}, \qquad 
  h=\frac{-\beta'_2+\beta'_4}{\beta'_1+\beta'_3}.
$$
Hence,
$$
	\beta'_2=f(g-h), \qquad    \beta'_4=f(g+h).
$$
Note that the function $f$ doesn't vanish and $g$, $h$ are well defined because $	(\beta_1')^2-(\beta_3')^2\ne 0 $.
Using (\ref{eq:3.1}) we get also 
$$
	\beta'_1=f(1+gh), \qquad \beta'_3=f(1-gh).
$$
Finally, we obtain that  if $(\beta_1')^2-(\beta_3')^2$ doesn't vanish in $J_0$, then 
\begin{equation} \label{eq:3.1-a}
    \beta'(u)=f(u)\big(1+g(u)h(u), g(u)-h(u), 1-g(u)h(u), g(u)+h(u) \big), \quad u \in J_0.
\end{equation} 

\vskip 2mm
If $(\beta_1')^2-(\beta_4')^2 \neq 0$ at a point $u_0 \in J$, then there exists $J_0 \subset J$ such that  the function
$(\beta_1')^2-(\beta_4')^2 \neq 0$ in $J_0$. In such case we can introduce in $J_0$ functions:
$$
	f=\frac{\beta'_1+\beta'_4}{2}, \qquad 	g=\frac{\beta'_2+\beta'_3}{\beta'_1+\beta'_4}, \qquad 
h=\frac{-\beta'_2+\beta'_3}{\beta'_1+\beta'_4}.
$$
Using (\ref{eq:3.1}) we get
$$
	\beta'_1=f(1+gh), \qquad \beta'_2=f(g-h), \qquad \beta'_3=f(g+h), \qquad    \beta'_4= f(1-gh).
$$
So, in this case  we obtain 
\begin{equation} \label{eq:3.1-b}
    \beta'(u)=f(u)\big(1+g(u)h(u), g(u)-h(u), g(u)+h(u), 1-g(u)h(u) \big), \quad u \in J_0.
\end{equation} 
Now, we substitute
$$
		\tilde{g} = \frac{1-h}{1+h}, \qquad \tilde{h} = \frac{1-g}{1+g}
$$
and get 
$$
		g-h = \frac{2 (\tilde{g} - \tilde{h})}{(1+\tilde{h}) (1+\tilde{g})}, \qquad g+h = \frac{2 (1 - \tilde{g}\tilde{h})}{(1+\tilde{h}) (1+\tilde{g})}, 
$$
$$
		1 + gh = \frac{2 (1 + \tilde{g}\tilde{h})}{(1+\tilde{h}) (1+\tilde{g})}, \qquad 1 - gh = \frac{2 (\tilde{g} + \tilde{h})}{(1+\tilde{h}) (1+\tilde{g})}.
$$
Setting 
$$
	\tilde{f}= \frac{2 f}{(1+\tilde{h}) (1+\tilde{g})}
$$
we get
\begin{equation} \notag
    \beta'(u)=\tilde{f}(u)\big(1+\tilde{g}(u)\tilde{h}(u), \tilde{g}(u)-\tilde{h}(u), 1-\tilde{g}(u)\tilde{h}(u), \tilde{g}(u)+\tilde{h}(u) \big).
\end{equation} 
Hence, omitting points at which $1+h=0$ or $1+g =0$ we can transform \eqref{eq:3.1-b} in \eqref{eq:3.1-a}.

Consequently, we have proved the following statement.

\vspace{0.2 cm}
\begin{prop}\label{T:3.2}   
Each null curve $\beta=\big( \beta_1,\beta_2,\beta_3,\beta_4\big)$  in $\R^4_2$ can be parametrized locally (in a neighbourhood of a given point) in such a way that

\vskip 1mm
(1) $\beta'(u)=f(u)\big( 1,\varepsilon_1,\varepsilon_2,\varepsilon_3\big), \quad \varepsilon_i=\pm 1, \; i = 1, 2, 3;$\\
or
\vskip 1mm
(2) $\beta'(u)=f(u)\big( 1+g(u)h(u),g(u)-h(u),1-g(u)h(u),g(u)+h(u) \big)$\\
for some smooth functions $f(u)$, $g(u)$, $h(u)$, such that $f(u)$ does not vanish.
\end{prop}

\vspace{0.2 cm}

In what follows we  use Propositions \ref{T:3.1}
and  \ref{T:3.2} to find a Weierstrass representation of a minimal Lorentzian surface in  $\R^4_2$ with two-dimensional first normal space.

\vskip 2mm
Let $S$ be a minimal surface in $\R^4_2$. According to Proposition \ref{T:3.1}, $S$  can be parametrized in isothermal parameters by
$$	S:  \Psi(u,v)=\beta(u+v)+\theta(u-v),$$
where $\beta$ and $\theta$ are null curves in  $\R^4_2$ such that $\langle\beta',\theta'\rangle \neq 0$.

\vskip 2mm
(i) Let us suppose that both  $\beta$ and $\theta$ are of  type (1) in Proposition \ref{T:3.2}, i.e. 
$$\beta' = f_1 \big( 1, \varepsilon_1, \varepsilon_2, \varepsilon_3\big), \qquad \theta' = f_2 \big( 1, \varepsilon_4, \varepsilon_5, \varepsilon_6\big), \quad \varepsilon_i=\pm 1, \; i = 1, \dots, 6.$$
We denote by $c_1$ and $c_2$ the constant vectors $\big( 1, \varepsilon_1, \varepsilon_2, \varepsilon_3\big)$ and $\big( 1, \varepsilon_4, \varepsilon_5, \varepsilon_6\big)$, respectively. 
Then
$$ \beta'(u+v) = f_1(u+v) \,c_1, \qquad \theta'(u-v) = f_2(u-v) \,c_2.$$
Hence, the tangent space of $S$ is spanned by the vector fields
$$\begin{array}{l}
\vspace{2mm}
\Psi_u = f_1(u+v)\, c_1 + f_2(u-v) \,c_2,\\
\vspace{2mm}
\Psi_v = f_1(u+v) \,c_1 - f_2(u-v)\, c_2
\end{array}
$$
and we get
$$\begin{array}{l}
\vspace{2mm}
\Psi_{uu} = f'_1(u+v)\,c_1 + f'_2(u-v) \,c_2,\\
\vspace{2mm}
\Psi_{uv} = f'_1(u+v) \,c_1 - f'_2(u-v) \,c_2,\\
\vspace{2mm}
\Psi_{vv} = f'_1(u+v)\, c_1 + f'_2(u-v) \,c_2.
\end{array}$$
The last three equalities show that $\sigma(\Psi_u, \Psi_u)=0$, $\sigma(\Psi_u, \Psi_v) = 0$, $\sigma(\Psi_v, \Psi_v)=0$, i.e. 
the surface $S$ has vanishing second fundamental tensor $\sigma$. Hence, $S$ is part of a plane.

\vskip 2mm
(ii) Let $\beta$ be of type (1) and $\theta$ be of  type (2) in Proposition \ref{T:3.2}, i.e.
$$\begin{array}{l}
\vspace{2mm}
\beta' = f_1 \,c,\\
\vspace{2mm}
\theta' =  f_2\big(1 + gh, g - h, 1 - gh, g + h \big),
\end{array}$$
where $c= \big( 1, \varepsilon_1, \varepsilon_2, \varepsilon_3\big), \varepsilon_i=\pm 1$, $\beta'$ depends on  $u+v$, $\theta'$ depends on  $u-v$.
Hence, the tangent space of $S$ is spanned by 
$$\begin{array}{l}
\vspace{2mm}
\Psi_u = f_1(u+v) \,c + \theta'(u-v),\\
\vspace{2mm}
\Psi_v = f_1(u+v) \,c - \theta'(u-v),
\end{array}
$$
and 
$$\begin{array}{l}
\vspace{2mm}
\Psi_{uu} = f'_1(u+v) \, c + \theta''(u-v),\\
\vspace{2mm}
\Psi_{uv} = f'_1(u+v) \,c - \theta''(u-v),\\
\vspace{2mm}
\Psi_{vv} = f'_1(u+v) \,c + \theta''(u-v).
\end{array}$$
The above equalities imply that $\sigma(\Psi_u, \Psi_u)= -\sigma(\Psi_u, \Psi_v) = \sigma(\Psi_v, \Psi_v)$, since $f'_1(u+v)\, c$ is a tangent vector field. Hence, the  first normal space of the surface is one-dimensional. 

\vskip 2mm
Consequently, when the first normal space of the surface is two-dimensional,  the null curves $\beta$ 
and $\theta$ are both of type (2) in Proposition \ref{T:3.2}. So, we have:

\begin{thm}\label{T:3.3}   
Up to change of the parameters any minimal Lorentzian surface  in $\R^4_2$ with two-dimensional first normal space has locally 
a Weierstrass type representation in isothermal parameters
$$
	\Psi(u,v)=\beta(u+v)+\theta(u-v)
$$
for some null curves $\beta$ and $\theta$ satisfying
$$\begin{array}{l}
\vspace{2mm}
    \beta'(t)=f_1(t)\big( 1+g_1(t)h_1(t), g_1(t)-h_1(t), 1-g_1(t)h_1(t), g_1(t)+h_1(t) \big);\\
\vspace{2mm}
    \theta'(s)=f_2(s)\big( 1+g_2(s)h_2(s), g_2(s)-h_2(s), 1-g_2(s)h_2(s), g_2(s)+h_2(s) \big), 
\end{array}$$
and $\langle\beta'(t),\theta'(s)\rangle \neq 0$.
\end{thm}


\vspace{.5 cm}

\setcounter{equation}{0}
\section{Minimal Lorentzian surfaces of general type and their canonical Weierstrass representation} \label{S:Canonical}

Let $S: \Psi=\Psi(u,v), \; (u, v) \in \mathcal{D}, \,  \mathcal{D} \subset \R^2$ be a minimal Lorentzian surface of general type. 
According to Theorem \ref{T:3.3}, $S$ is parametrized by 
\begin{equation} \label{Weierstrass}
\Psi(u,v)=\beta(u+v)+\theta(u-v),
\end{equation}
where $\beta$ and $\theta$ have the form
\begin{equation} \label{Weierstrass-deriv}
	\begin{array}{l}
	\vspace{2mm}
		\beta'(t)=f_1(t)\big( 1+g_1(t)h_1(t),g_1(t)-h_1(t),1-g_1(t)h_1(t),g_1(t)+h_1(t) \big), \\
		\vspace{2mm}
		\theta'(s)=f_2(s)\big( 1+g_2(s)h_2(s),g_2(s)-h_2(s),1-g_2(s)h_2(s),g_2(s)+h_2(s) \big).
	\end{array} 
\end{equation}
Then the tangent space is spanned by
\begin{equation} \label{E:first_der}
\begin{array}{l}
\vspace{2mm}
\Psi_u = \beta'(u+v) + \theta'(u-v),\\
\vspace{2mm}
\Psi_v = \beta'(u+v) - \theta'(u-v),
\end{array}
\end{equation}
and the coefficients $E$, $F$, and $G$ of the first fundamental form of $S$ are
$$
	\begin{array}{l}
	\vspace{2mm}
	E=4f_1(u+v)f_2(u-v)\big( g_1(u+v)-g_2(u-v)\big) \big( h_1(u+v)-h_2(u-v) \big), \\
	\vspace{2mm}
	F=0, \\
	\vspace{2mm}
	G=-E.
	\end{array}
$$

We consider the orthonormal tangent frame field  $x=\ds{\frac{\Psi_u}{\sqrt{E}}}$, $y=\ds{\frac{\Psi_v}{\sqrt{-G}}}$.
Note that (\ref{eq:2.1}) and (\ref{eq:2.2}) imply
$$
	K^2-\varkappa^2=(a^2-b^2)^2+(c^2-d^2)^2+2(a^2-b^2)(c^2-d^2)-4(ac-bd)^2
$$
or equivalently
\begin{equation} \label{eq:4.1}
	\begin{array}{rl}
	\vspace{2mm}
	 K^2-\varkappa^2= & \langle \sigma(x,x),\sigma(x,x) \rangle^2+\langle \sigma(x,y),\sigma(x,y) \rangle^2  \\
	\vspace{2mm}
	 & +2\langle \sigma(x,x),\sigma(x,x) \rangle\langle \sigma(x,y),\sigma(x,y) \rangle-4\langle \sigma(x,x),\sigma(x,y) \rangle^2. 
		\end{array}
\end{equation}

Using \eqref{E:first_der} we get
\begin{equation} 
\begin{array}{l}  \label{Eq:4.5}
	\vspace{2mm}
	\sigma(x,x) = \ds{\frac{1}{E} \left(\beta'' + \theta'' - \frac{\langle \beta'', \theta' \rangle}{\langle \beta', \theta' \rangle} \,\beta' - 
	\frac{\langle \beta', \theta'' \rangle}{\langle \beta', \theta' \rangle} \,\theta' \right)};\\
	\vspace{2mm}
		\sigma(x,y) = \ds{\frac{1}{E} \left(\beta'' - \theta'' - \frac{\langle \beta'', \theta' \rangle}{\langle \beta', \theta' \rangle}\, \beta' + 
	\frac{\langle \beta', \theta'' \rangle}{\langle \beta', \theta' \rangle}\,\theta' \right)}.
	\end{array}
\end{equation}
On the other hand, from \eqref{Weierstrass-deriv} it follows that 
\begin{equation} \label{Eq:4.6}
\begin{array}{l}
	\vspace{2mm}
	\langle \beta', \beta'' \rangle = 0;\\
	\vspace{2mm}
		\langle \theta', \theta'' \rangle = 0;\\
	\vspace{2mm}
		\langle \beta'', \beta'' \rangle = -4 (f_1)^2 g_1'h_1';\\
	\vspace{2mm}
		\langle \theta'', \theta'' \rangle = -4 (f_2)^2 g_2'h_2';\\
		\vspace{2mm}
		\langle \theta'', \beta' \rangle = 2 f_1 f_2' (g_2-g_1)(h_2-h_1) + 2 f_1 f_2 ((g_2-g_1)h_2' + (h_2-h_1)g_2');\\
	\vspace{2mm}
		\langle \beta'', \theta' \rangle = 2 f_2 f_1' (g_2-g_1)(h_2-h_1) + 2 f_1 f_2 ((h_1-h_2)g_1' + (g_1-g_2)h_1').\\
	\end{array}
\end{equation}

Now, using \eqref{eq:4.1}, \eqref{Eq:4.5} and \eqref{Eq:4.6}, by a straightforward calculation  we obtain
\begin{equation} \label{eq:4.3}
	K^2\!-\!\varkappa^2=\frac{g_1'(u\!+\!v)g_2'(u\!-\!v)h_1'(u\!+\!v)h_2'(u\!-\!v)}
	                     {(f_1(u\!+\!v))^2 (f_2(u\!-\!v))^2\big(g_1(u\!+\!v)-g_2(u\!-\!v)\big)^4\big(h_1(u\!+\!v)-h_2(u\!-\!v)\big)^4}.
\end{equation}

Note that the condition $K^2-\varkappa^2 \neq 0$ and (\ref{eq:4.3}) imply that none of the 
functions $ g_1'$, $g_2'$, $h_1'$, $h_2' $ can vanish. Moreover, since $K^2-\varkappa^2>0$ for  minimal Lorentzian surfaces of general type, then $g_1'g_2' h_1'h_2' >0$. Hence, as a consequence of  Theorem \ref{T:3.3} we obtain the  Weierstrass representation of minimal Lorentzian surfaces of general type given below.

	\begin{thm}\label{T:3.3-general type}   
Up to change of the parameters any minimal Lorentzian surface of general type in $\R^4_2$ has locally 
a Weierstrass type representation in isothermal parameters
$$
	\Psi(u,v)=\beta(u+v)+\theta(u-v)
$$
for some null curves $\beta$ and $\theta$ satisfying
$$\begin{array}{l}
\vspace{2mm}
    \beta'(t)=f_1(t)\big( 1+g_1(t)h_1(t), g_1(t)-h_1(t), 1-g_1(t)h_1(t), g_1(t)+h_1(t) \big);\\
\vspace{2mm}
    \theta'(s)=f_2(s)\big( 1+g_2(s)h_2(s), g_2(s)-h_2(s), 1-g_2(s)h_2(s), g_2(s)+h_2(s) \big), 
\end{array}$$
and $f_1 f_2 (g_1-g_2)(h_1-h_2) \neq 0$, $g_1'g_2' h_1'h_2' >0$ everywhere.
\end{thm}

\vskip 3mm
Further, we shall obtain a condition for the surface to be parametrized by canonical parameters. A necessary condition for a canonical parametrization of a minimal surface of general type is $\sigma(\Psi_{u},\Psi_{u})$ and $\sigma(\Psi_{u},\Psi_{v})$ to be orthogonal  \cite{A-M}. It follows from \eqref{Eq:4.5} and \eqref{Eq:4.6} that
$$
	\begin{array}{l}
	\vspace{2mm}
		\langle \sigma(\Psi_{u},\Psi_{u}),\sigma(\Psi_{u},\Psi_{v}) \rangle \\
		\vspace{2mm}
		 =\langle\beta'',\beta''\rangle(u+v)-\langle\theta'',\theta''\rangle(u-v)\\ 
		\vspace{2mm}
		=-4(f_1(u+v))^2g_1'(u+v)h_1'(u+v)+4(f_2(u-v))^2g_2'(u-v)h_2'(u-v).
	\end{array}
$$
Hence, for a canonical parametrization we must have 
$$
	(f_1(u+v))^2g_1'(u+v)h_1'(u+v) = (f_2(u-v))^2g_2'(u-v)h_2'(u-v)=const \neq 0.
$$ 
Suppose that 
$$
	f_1=\frac1{4\sqrt{|g_1'h_1'|}}, \qquad f_2=\frac1{4\sqrt{|g_2'h_2'|}}.
$$
Then,  using   \eqref{eq:2.6} and \eqref{eq:4.3} we see that
$$
	\left|\mu^2-\nu^2\right|=16
	\frac{g_1'(u+v)h_1'(u+v)g_2'(u-v)h_2'(u-v)}
	{\big( g_1(u+v)-g_2(u-v)\big)^2\big( h_1(u+v)-h_2(u-v)\big)^2}.
$$
So,  \eqref{eq:2.4} also holds and hence, the surface is parametrized by canonical parameters.
Consequently, we have proved:

\begin{prop}\label{T:4.2}  
Let $\beta$ and $\theta$ be null curves in $\R^4_2$  satisfying
\begin{equation} \label{canon-deriv}
	\begin{array}{l}
    \ds \beta'(t)=\frac1{4\sqrt{|g_1'(t)h_1'(t)|}}\big( 1\!+\!g_1(t)h_1(t),g_1(t)\!-\!h_1(t),1\!-\!g_1(t)h_1(t),g_1(t)\!+\!h_1(t) \big) \\
    \ds \theta'(s)=\frac1{4\sqrt{|g_2'(s)h_2'(s)|}}\big( 1\!+\!g_2(s)h_2(s),g_2(s)\!-\!h_2(s),1\!-\!g_2(s)h_2(s),g_2(s)\!+\!h_2(s) \big)  
	\end{array}
\end{equation}
with $(g_1-g_2)(h_1-h_2) \neq 0$, $g_1'g_2' h_1'h_2' >0$. Then 
\begin{equation} \label{canon}
	\Psi(u,v)=\beta(u+v)+\theta(u-v)
\end{equation}
is a minimal Lorentzian surface of general type in $\R^4_2$ parametrized by canonical parameters.
\end{prop}

\vskip 3mm
Now we want to see whether it is possible to change the parameters of a minimal surface 
given in a Weierstrass form so that the obtained representation be canonical. Consider a minimal surface given by (\ref{Weierstrass}), such that the null
curves $\beta$ and $\theta$ satisfy (\ref{Weierstrass-deriv}) and $f_1 f_2 (g_1-g_2)(h_1-h_2) \neq 0$, $g_1'g_2' h_1'h_2' >0$. We look
for a change of the parameters 
$$
	u=p(\bar u+\bar v)+q(\bar u-\bar v) \ , \qquad
	v=p(\bar u+\bar v)-q(\bar u-\bar v) 
$$
for some functions $p$, $q$ of one variable, such that  the function
$$
	\overline\Psi(\bar u, \bar v)=\Psi\big( p(\bar u+\bar v)+q(\bar u-\bar v),p(\bar u+\bar v)-q(\bar u-\bar v) \big)
$$
has the form 
\begin{equation} \label{bar-canon}
	\overline\Psi(\bar u,\bar v)=\bar\beta(\bar u+\bar v)+\bar\theta(\bar u-\bar v)
\end{equation}
for some functions $\bar\beta$ and $\bar\theta$ satisfying
\begin{equation} \label{bar-canon-deriv}
	\begin{array}{l}
    \ds \bar\beta'(\bar u)\!=\!
		    \frac1{4\sqrt{|\bar g_1'(\!\bar u\!)\bar h_1'(\!\bar u\!)|}}\big(\!1\!+\!\bar g_1(\bar u)\bar h_1(\bar u),
				\bar g_1(\bar u)\!-\!\bar h_1(\bar u),1\!-\!\bar g_1(\bar u)\bar h_1(\bar u),\bar g_1(\bar u)\!+\!\bar h_1(\bar u)\!\big), \\
    \ds \bar \theta'(\bar v)\!=\!\frac1{4\sqrt{|\bar g_2'(\bar v)\bar h_2'(\bar v)|}}\big(\!1\!+\!\bar g_2(\bar v)\bar h_2(\bar v),
		    \bar g_2(\bar v)\!-\!\bar h_2(\bar v),1\!-\!\bar g_2(\bar v)\bar h_2(\bar v),\bar g_2(\bar v)\!+\!\bar h_2(\bar v)\!\big),  
	\end{array}
\end{equation}
and as a consequence of Proposition \ref{T:4.2} this
will be a representation of the surface in canonical parameters. It 
can be seen that
$$
	\begin{array}{l}
	\vspace{2mm}
		\overline\Psi_{\bar u}=(p'+q')\Psi_u+(p'-q')\Psi_v \\
		\vspace{2mm}
		\overline\Psi_{\bar v}=(p'-q')\Psi_u+(p'+q')\Psi_v
	\end{array}
$$
imply
$$
	\bar\beta'(z)=2\beta'(p(z))p'(z), \qquad \bar\theta'(z)=2\theta'(q(z))q'(z).
$$
Hence, we derive
\begin{equation} \label{eq:4.4}
	\begin{array}{l}
		\ds (p'(z))^2=\frac{1}{8 |f_1(p(z))|\sqrt{|g_1'(p(z))h_1'(p(z))|}}, 
					\quad \bar g_1(z)=g_1(p(z)), \;\; \bar h_1(z)=h_1(p(z)); 
					\vspace{3mm} \\
		\ds (q'(z))^2=\frac{1}{8 |f_2(q(z))|\sqrt{|g_2'(q(z))h_2'(q(z))|}}, 
					\quad \bar g_2(z)=g_2(q(z)), \;\; \bar h_2(z)=h_2(q(z)).
	\end{array}
\end{equation}
Consequently, we may determine not only the change of the parameters, i.e. the functions $p(z)$, $q(z)$, but also 
a canonical representation of the given surface. Finally, we proved:

\begin{thm}\label{T:4.3}  
Let $S$ be a minimal Lorentzian surface of general type  in $\R^4_2$ parametrized by (\ref{Weierstrass}), 
such that the null curves $\beta$ and $\theta$ satisfy (\ref{Weierstrass-deriv}),
where  $f_1 f_2 (g_1-g_2)(h_1-h_2) \neq 0$, $g_1'g_2' h_1'h_2' >0$. Then, $S$ has a canonical representation
(\ref{bar-canon}) where $\bar\beta$ and $\bar\theta$ are null curves satisfying (\ref{bar-canon-deriv}).
This canonical representation is determined by  equations (\ref{eq:4.4}).
\end{thm}

\vspace{0.2cm}
We can also  find the Gauss curvature $K$ and the curvature of the normal connection  $\varkappa$ 
of a minimal Lorentzian surface of general type  in $\R^4_2$ given in canonical
parameters. Namely, for the surface (\ref{canon}), where the null curves
$\beta$ and $\theta$ are given by (\ref{canon-deriv}), we obtain:

\begin{equation} \label{eq:4.5}
	\begin{array}{l}
	\vspace{3mm}
		\ds K=-8 \frac{\sqrt{g_1'h_1'g_2'h_2'} 
			\big( (h_1-h_2)^2g_1'g_2'+(g_1-g_2)^2h_1'h_2' \big)}
			{(g_1-g_2)^3(h_1-h_2)^3};\\  
			\vspace{2mm}
		\ds \varkappa=-8 \sqrt{g_1'h_1'g_2'h_2'}
			\left|\frac{ (h_1-h_2)^2g_1'g_2'-(g_1-g_2)^2h_1'h_2' }
			{(g_1-g_2)^3(h_1-h_2)^3} \right|,  
		\end{array}	
\end{equation}
where  $g_1,\, h_1$ depend on $u+v$ and
$g_2,\, h_2$ depend on $u-v$. Hence, using Theorem B we obtain

\begin{thm}\label{T:4.4}  
The system of partial differential equations  (\ref{eq:2.5}) has solutions of type
(\ref{eq:4.5})
for arbitrary smooth functions $g_1,g_2,h_1,h_2$, such that $g_1'g_2' h_1'h_2' >0$ 
and the differences $g_1-g_2$ and $h_1-h_2$ never vanish.
\end{thm}


\vspace{.5 cm}

\setcounter{equation}{0}
\section{Examples} 

In this section, using Theorem {\ref{T:3.3}} we give some examples of 
minimal Lorentzian surfaces   in $\R^4_2$. Moreover,
according to Theorem {\ref{T:3.3-general type}} they are of general type and
according to Proposition {\ref{T:4.2}} they are in canonical
parametrization. 

\begin{example}\label{Ex:5.1}
Let $g_1(t)=t+1$, $h_1(t)=t$, $g_2(s)=s$, $h_2(s)=s+1$.  Using  (\ref{canon-deriv}) we have
\begin{equation} \notag
	\begin{array}{l}
	\vspace{2mm}
    \ds \beta'(t)=\frac1{4}\big( 1+ t (t+1); 1; 1 - t (t+1); 2t +1 \big); \\
	\vspace{2mm}
    \ds \theta'(s)=\frac1{4}\big( 1+ s (s+1); - 1; 1 - s (s+1); 2s +1\big).  
	\end{array}
\end{equation}
Then up to constant vectors the null curves $\beta$ and $\theta$ are given by
\begin{equation} \notag
	\begin{array}{l}
	\vspace{2mm}
    \ds \beta(t)=\frac1{24}\big(2t^3+3t^2+6t; 6t; -2t^3-3t^2+6t; 6t^2 +6t\big); \\
	\vspace{2mm}
    \ds \theta(s)=\frac1{24}\big( 2s^3+3s^2+6s; -6s; -2s^3-3s^2+6s; 6s^2 +6s\big).  
	\end{array}
\end{equation}
Applying  (\ref{canon}) we obtain the surface 
$$
	\Psi(u,v)=\left\{
	                 \begin{array}{l}
									\vspace{2mm}
	                    \ds\frac1{12}\big(2u^3+3u^2+3v^2+6u(v^2+1)\big)  \smallskip \\
											\vspace{2mm}
	                    \ds \frac v2 \smallskip \\
											\vspace{2mm}
	                    \ds \frac1{12}\big(-2u^3-3u^2-3v^2-6u(v^2-1)\big)\smallskip \\
											\vspace{2mm}
	                    \ds \frac12\big(u+u^2+v^2\big) \ .
	                 \end{array}\right.
$$
In this example, $(g_1(u+v)-g_2(u-v)) (h_1(u+v) - h_2(u-v)) = (2v+1) (2v-1)$, so we assume that $v \neq \pm \frac{1}{2}$.

According to (\ref{eq:4.5}) the Gauss curvature $K$ and the curvature  of the normal connection $\varkappa$ are expressed by:
$$
	K=\frac{16(1+4v^2)}{(1-4v^2)^3} \ ; \qquad \varkappa=-\left|\frac{64v}{(1-4v^2)^3}\right| \ .
$$
\end{example}

\vskip 3mm

\begin{example}\label{Ex:5.2}
Let $g_1(t)=t$, $h_1(t)=t$, $g_2(s)=e^s$, $h_2(s)=-e^{-s}$.  Using  (\ref{canon-deriv}) and  (\ref{canon}) we obtain the surface  
$$
	\Psi(u,v)=\left\{
	                 \begin{array}{l}
									\vspace{2mm}
	                    \ds\frac{1}{12}(u+v)^3 + \frac{1}{4}(u+v) \smallskip \\
											\vspace{2mm}
	                    \ds \frac 12 \sinh(u-v) \smallskip \\
											\vspace{2mm}
	                    \ds -\frac{1}{12}(u+v)^3 + \frac{1}{4}(3u-v) \smallskip \\
											\vspace{2mm}
	                    \ds \frac{1}{2}\cosh(u-v)+\frac{1}{4}(u+v)^2.
	                 \end{array}\right.
$$
In this example, $(g_1(u+v)-g_2(u-v)) (h_1(u+v) - h_2(u-v)) = (u+v - e^{u-v}) (u+v + e^{-(u-v)})$, so we assume that the parameters 
$u$ and $v$ satisfy $u+v \neq e^{u-v}$, $u+v \neq -e^{-(u-v)}$, or equivalently $(u+v)^2 - 1 \neq 2(u+v) \sinh (u-v)$.

The Gauss curvature $K$ and the curvature of the normal connection $\varkappa$  are:
$$
	K= -16\frac{\cosh(u-v)\big((u+v)^2+1\big)}
	         {\big((u+v)^2 - 1 -2(u+v) \sinh (u-v)\big)^3} \ , 
$$
$$
	\varkappa=-16\left|\frac{\big((u+v)^2-1\big) \sinh(u-v) + 2(u+v)}
                 {\big((u+v)^2 - 1 -2(u+v) \sinh (u-v)\big)^3}\right| \ .
$$
\end{example}

Note, that according to Theorem \ref{T:4.4}, the functions $K$ and $\varkappa$ from Example \ref{Ex:5.1}
and Example \ref{Ex:5.2} give solutions to the system of partial differential equations (\ref{eq:2.5}).

\vskip 5mm \textbf{Acknowledgments:}
The  second author is partially supported by the National Science Fund, Ministry of Education and Science of Bulgaria under contract DN 12/2.

\end{document}